\documentclass[a4paper,11pt]{article}
\pdfoutput=1 % if your are submitting a pdflatex (i.e. if you have
\usepackage{jheppub} % for details on the use of the package, please
\usepackage[T1]{fontenc} % if needed

\usepackage{pdfsync}
\usepackage{mathptmx}
\usepackage[utf8]{inputenc}
\usepackage[T1]{fontenc}
\usepackage{slashed}
\usepackage{times}
\usepackage{amssymb,amsfonts,amsmath,amsthm}
\usepackage{dsfont,bbm}
\usepackage{dcolumn}
\usepackage{epsf}
\usepackage{graphicx}
\usepackage[caption=false]{subfig}
\usepackage{dsfont}
\usepackage{simplewick}
\usepackage{bm}
\usepackage{eucal}
\usepackage{slashed}
\usepackage[active]{srcltx}
\usepackage[usenames]{color}
\title{\boldmath Universal inverse Radon transforms: Inhomogeneity, angular restrictions and boundary
}

\author[a]{I.~V.~Anikin}
\affiliation[a]{Bogoliubov Laboratory of Theoretical Physics JINR, 141980 Dubna, Russia}

\emailAdd{anikin@theor.jinr.ru}

\abstract{
An alternative method to invert the Radon transforms without the use of Courant-Hilbert's identities
has been proposed and developed independently from the space dimension.
For the universal representation of inverse Radon transform, we study
the consequences of inhomogeneity of outset function without the restrictions
on the angular Radon coordinates.
We show that this inhomogeneity yields a natural evidence for
the presence of the extra contributions in the case of the full angular region.
In addition, if the outset function is well-localized in the space,
we demonstrate that the corresponding boundary conditions and the angular restrictions
should be applied for both the direct and inverse Radon transforms.
Besides, we relate the angular restrictions on the Radon variable
to the boundary exclusion of outset function and its Radon image.
}

\begin{document}
\maketitle
\flushbottom
%%%%%%%%%%%%%%%%%%%%%%%%%%%%%%%%%%%%%%%%%%

%%%%%%%%%%%%%%%%
\section{Introduction}
\label{Intro}
%%%%%%%%%%%%%%%%%%%%%%%%%%%%%%%%%%%%

Nowadays, there is no need to explain the importance of computerized tomography (CT) which influences on the different fields
and gives a possibility to investigate the internal composite structure of a object without breaking it.
The mathematical foundation of CT is associated with both
the direct and inverse Radon transforms \cite{Deans}. The inverse Radon transforms are being used on the stage
of visualization
\footnote{In hadron physics, the inverse Radon transform is used in processing data from the
generalized parton distributions  \cite{Belitsky:2000vk}.}. 
Meanwhile, the inversion of Radon transforms meets the problems related to the ill-posedness,
see for example \cite{Anikin:2019oes, Anikin:2024vto}.
The explicit expression of inversion depends on the (even or odd) dimension of space where the outset function, which
has to be reconstructed, is defined \cite{Deans}. On the other hand,
there are no the well-defined methods to derive the universal form of Radon inversion
that is suitable for both even and odd dimension simultaneously.
The standard methods of inversion for even and odd dimension have been based on
the use of Courant-Hilbert's identities which have the different forms depending on the space dimension \cite{Courant-Hilbert}.
In Courant-Hilbert's identities, the angular integration has been always performed over
the full region of variations, {\it i.e.} in the full interval $(0,\, 2\pi)$.
Besides, in the standard methods, the choice of the angular interval of integration has been usually dictated only by the
corresponding normalization constants and it has no much a deep physical (or/and mathematical) meaning.

In the present paper, we derive the universal inverse Radon transforms
without the use of Courant-Hilbert's identities.
The proposed method is backed by the
regularization within the generalized function theory \cite{Gelfand:1964, GGV}.

We demonstrate that the universal inverse Radon transform involves two
essential terms $f_S$ and $f_A$ which do contribute, independently from the space dimension,
even in the full angular region, $(0, 2\pi)$, owing to the inhomogeneity property of outset function.

We also find that if the outset function is well-localized in the space,
we are forced to deal with the restrictions which have been imposed on the Radon angular dependence.
We implement the detail analysis of how the restricted angular dependence of Radon transforms
appears as a consequence of the finite support of outset functions localized in the space domain.
These angular restrictions give the other evidences supporting the existence of two essential terms in the
universal inversion of Radon transformations.
It is worth to notice that the two mentioned terms, $f_S$ and $f_A$, lead to the complexity of Radon inversion.
In its turns, this complexity opens a possibility to extend and to improve the Tikhonov regularization
needed for the different practical applications \cite{Anikin:2024vto}.

In the paper, we explore the other key moments which are related to
the connections between the angular restrictions, the presence of surface terms in the corresponding integrations
and the boundary exclusion of outset function (and, therefore, of the Radon images).

%%%%%%%%%%%%%%%%
\section{Basis of universal inverse Radon transforms}
\label{IRT-1}
%%%%%%%%%%%%%%%%%%%%%%%%%%%%%%%%%%%%

As well-known, in order to restore the needed information on the given internal structure
described by the outset function $f(\vec{\bf x})$ ($\vec{\bf x}\in \mathbb{R}^n$ for $\forall n$), one has to use
the inverse Radon transform
that expresses the outset function $f(\vec{\bf x})$ through the Radon image.
Namely, the direct Radon transform of $f(\vec{\bf x})$ is defined by 
\begin{eqnarray}
\label{F-t-4-dir}
G\equiv{\cal R}[f](\tau, \varphi, \theta_i) \quad \Longrightarrow \quad
\mathcal{R}[f](\tau,\varphi, \theta_i)=
\int_{-\infty}^{+\infty} d^n \vec{\bf x} \, f(\vec{\bf x})
\delta\left( \tau - \langle\vec{\bf n}_{\varphi, \theta_i}, \vec{\bf x}\rangle\right),
\end{eqnarray}
where the $n$-dimensional vector
$\vec{\bf n}_{\varphi, \theta_i}$ corresponds to the unit vector pointing along the radial Radon coordinate.
While the inverse Radon transform can be expressed as
\begin{eqnarray}
\label{Intro-1}
f ={\cal R}^{-1} G \quad \Longrightarrow \quad
f(\vec{\bf x}) =
\int d\mu_n(\eta, \varphi, \theta_i)\,
{\cal R}[f]\big(\eta + \langle \vec{\bf n}_{\varphi, \theta_i}, \vec{\bf x}\rangle, \varphi, \theta_i\big),
\end{eqnarray}
where the integration measure depends on the radial $\eta$ and the angular $\varphi, \theta_i$ Radon coordinates.
The explicit forms of integration measure have been specified below (see, for example, (\ref{Inv-F-t-2-2}), (\ref{Inv-F-t-2-3}) 
and (\ref{Inv-F-t-2-4-S2}), (\ref{Inv-F-t-2-4-A2})).
Also, it can involve the corresponding weight operator. 
Notice that, for the practical aims, our consideration deals only with the well-localized outset functions $f(\vec{\bf x})$
defined on the finite support. As a result, the radial Radon coordinate, $\tau$, should be finite-restricted too. 

In the standard approach to the inversion of Radon transforms \cite{Deans}, the different forms of $d\mu_n$ are specified by the space dimension which is either even or odd. Moreover, it is crucially important that the angular coordinate of the standard
inversion varies in the full region without any restrictions owing to the use of Courant-Hilbert's identities \cite{Courant-Hilbert}.

In this section, based on the methods of generalized function theory,
we derive the universal representation for the inverse Radon transform that is valid for an arbitrary
(for both even and odd) dimension of space. The universal representation stems from the Fourier slice theorem asserting that
{\it the Fourier image of the outset function $f(\vec{\bf x})$ relates to the direct Radon image of the same
function $f(\vec{\bf x})$ through the one-dimensional Fourier transformation with respect to the radial coordinate}:
\begin{eqnarray}
\label{F-t-3-dir}
\boxed{
\mathcal{F}[f](\lambda, \varphi, \theta_i) =\int_{-\infty}^{+\infty}(d\tau) \,e^{-i\lambda\tau}\,
\mathcal{R}[f](\tau,\varphi, \theta_i),
}
\end{eqnarray}
where the Fourier image of $f(\vec{\bf x})$ is given by
($\vec{\bf q}\equiv \lambda\, \vec{\bf n}_{\varphi, \theta_i}$ with $|\vec{\bf n}_{\varphi, \theta_i}|=1$)
\begin{eqnarray}
\label{F-t-1-dir}
\mathcal{F}[f](\vec{\bf q})= \int_{-\infty}^{+\infty} d^n \vec{\bf x} \,
e^{-i\langle\vec{\bf q},\vec{\bf x}\rangle} \,f(\vec{\bf x}).
\end{eqnarray}

In (\ref{F-t-3-dir}), the integration measure $(d\tau)$ includes the
corresponding normalization factor which is not written explicitly unless it leads to misunderstanding.

The principal inference of (\ref{F-t-3-dir}) is that the angular dependence of Fourier and Radon images are
coinciding. This fact is going to be used for our further derivations.

Using the inverse Fourier transform together with (\ref{F-t-3-dir}), we can write down that
\begin{eqnarray}
\label{Inv-F-t-2}
&&f(\vec{\bf x}) =
\int_{-\infty}^{+\infty} d^n \vec{\bf q}\,  e^{+i\langle\vec{\bf q},\vec{\bf x}\rangle} \,\mathcal{F}[f](\vec{\bf q})
\Big|_{\vec{\bf q}=\lambda\vec{\bf n}_{\varphi, \theta_i}}=
\\
\label{Inv-F-t-2-2}
&&
\int_{0}^{+\infty} d\lambda \lambda^{n-1} \int_{\text{f. r.}} d^{n-1}\Theta(\varphi, \theta_i)\,
e^{+i\lambda\langle \vec{\bf n}_{\varphi, \theta_i}, \vec{\bf x}\rangle}\,
\mathcal{F}[f](\lambda, \varphi, \theta_i).
\end{eqnarray}
Hence, the Fourier slice theorem, see (\ref{F-t-3-dir}), yields 
\footnote{$\epsilon$ as a subscript of $f$ denotes the $\epsilon$-regularization that
should be used in (\ref{Inv-F-t-2-3}).}
\begin{eqnarray}
\label{Inv-F-t-2-3}
&&
\hspace{-0.5cm}f_\epsilon(\vec{\bf x}) =
\int_{\text{f. r.}} d^{n-1}\Theta(\varphi, \theta_i)\,
\int_{-\infty}^{+\infty} (d\eta)\,
\mathcal{R}[f](\eta + \langle \vec{\bf n}_{\varphi, \theta_i}, \vec{\bf x}\rangle, \varphi, \theta_i)
\int_{0}^{+\infty} d\lambda \lambda^{n-1} \, e^{-i\lambda\eta}\Big|_{\epsilon\text{-reg.}},
\end{eqnarray}
where $``\epsilon\text{-reg.}"$ denotes the necessary regularization (see below) and
$``\text{f. r.}"$ signals that the angular integration measure covers the full regions of variations.
Generally speaking, it is already clear that (\ref{Inv-F-t-2-3}) establishes a basis for the universal inversion of Radon transforms.

As mentioned the representation of (\ref{Inv-F-t-2-3}) demands the $\epsilon$-regularization of the integration over
$\lambda$-variable. 
To this goal, we
make a replacement as $\eta\to\eta - i\epsilon$ which provides the analytical regularization
known from the distribution (generalized) function theory \cite{Gelfand:1964}.
In this case, the $\lambda$-integration reads
\begin{eqnarray}
\label{int-lam}
&&i^{\,2-n}\int_{0}^{+\infty} d\lambda \lambda^{n-1} \, e^{-i\lambda(\eta-i\epsilon)}=
i\,\frac{\partial^{n-1}}{\partial\eta^{n-1}} \int_{0}^{+\infty} d\lambda \, e^{-i\lambda(\eta-i\epsilon)}
\nonumber\\
&&
%\equiv 2\pi\frac{\partial^{n-1}}{\partial\eta^{n-1}} \delta_-(\eta)
=
 (-)^{n-1} (n-1)! \frac{\mathcal{P}}{\eta^n} + i\pi \frac{\partial^{n-1}}{\partial\eta^{n-1}}\delta(\eta),
\end{eqnarray}
where $\lambda$ as a pre-exponential factor has been traded for the derivative over $\eta$ acting on the exponential function.
Notice that, with the help of the $\epsilon$-regularization and (\ref{int-lam}), the inverse Radon representation of (\ref{Inv-F-t-2-3})
becomes to be well-defined in a sense of the regular (principle value) and singular (delta-function) generalized function/functional \cite{Gelfand:1964}.

It now remains to insert (\ref{int-lam}) into (\ref{Inv-F-t-2-3}) to obtain the universal inverse Radon transforms, we have
\begin{eqnarray}
\label{Inv-F-t-2-4}
\boxed{
f_\epsilon(\vec{\bf x}) = f_S(\vec{\bf x}) + f_A(\vec{\bf x}),
}
\end{eqnarray}
where
\begin{eqnarray}
\label{Inv-F-t-2-4-S}
&&
\hspace{-0.5cm}
f_S(\vec{\bf x}) =
i^{\,n-2}(-)^{n-1} (n-1)!\, \int_{\text{f.r.}} d^{n-1}\Theta(\varphi, \theta_i)\,
\int_{-\infty}^{+\infty} (d\eta)\, \frac{\mathcal{P}}{\eta^n}\,
\mathcal{R}[f](\eta + \langle \vec{\bf n}_{\varphi, \theta_i}, \vec{\bf x}\rangle, \varphi, \theta_i)
\end{eqnarray}
and
\begin{eqnarray}
\label{Inv-F-t-2-4-A}
&&
\hspace{-0.8cm}
f_A(\vec{\bf x})= (-)^{n-1}
i^{\,n-1}\, \pi \int_{\text{f.r.}} d^{n-1}\Theta(\varphi, \theta_i)\,
\int_{-\infty}^{+\infty} (d\eta)\,\delta(\eta)\,
\frac{\partial^{n-1}}{\partial\eta^{n-1}}\,\mathcal{R}[f](\eta + \langle \vec{\bf n}_{\varphi, \theta_i}, \vec{\bf x}\rangle, \varphi, \theta_i).
\end{eqnarray}
In (\ref{Inv-F-t-2-4-A}), we suppose $\mathcal{R}[f](\eta; ...)$ to be a restricted function of $\eta$, see below.

If there are no the angular restrictions, {\it i.e.} we deal with the full regions of angular integrations,
and the outset function $f(\vec{\bf x})$ is a homogeneous function with the trivial holonomy (that is, no holonomy at all),
the term $f_S$ gives the contribution only to the even dimension of space while the term $f_A$ contributes only to the
case of the odd dimension.
Indeed, to demonstrate this statement it is enough (a) to make the simultaneous replacements: $\eta\to -\eta$,
$\vec{\bf n}_{\varphi, \theta_i}\to - \vec{\bf n}_{\varphi, \theta_i}$ in (\ref{Inv-F-t-2-4-S}) and (\ref{Inv-F-t-2-4-A})
and (b) to take into account the corresponding symmetry properties of $\eta$-dependent coefficient functions.
This our observation reproduces the standard results for the inversion of Radon transforms,
see \cite{Deans}.

Notice that, within our approach based on (\ref{Inv-F-t-2-4}), the complexity of $f_A$ appears by the natural way
without an uncertainty. It is due to the use of Cauchy's theorem, see \cite{Anikin:2024vto}.
We stress that in the standard methods used up to now, the complexity of the inverse Radon transform for the odd
(or even depending on the precise normalization) dimension hidden in the corresponding irrelevant normalization
and it does not matter much.
It is worth to mention that the natural evidences for the complexity of  Radon transforms have been investigated 
in \cite{Anikin-CRT} owing to the introduction of the Wigner-like hybrid outset functions .

From (\ref{Inv-F-t-2-4-S}) and (\ref{Inv-F-t-2-4-A}), one can see that if the angular dependence
has been restricted for some reasons:
 \begin{eqnarray}
 \label{rest-ang}
 \boxed{
 \int_{\text{f.r.}} d^{n-1}\Theta(\varphi, \theta_i) \Longrightarrow \int_{\text{rest.r}} d^{n-1}\Theta(\varphi, \theta_i),
 }
 \end{eqnarray}
together with the broken (thanks for the inhomogeneity) symmetry of outset function,
{\it both terms $f_S$ and $f_A$ do contribute in the final inversion of Radon transforms independently from the space dimension}.
We prove this inference in the next sections.

%%%%%%%%%%%%%%%%
\section{The inhomogeneity of outset function }
\label{Non-homog}
%%%%%%%%%%%%%%%%%%%%%%%%%%%%%%%%%%%%

In this section, we study the role of the inhomogeneity of outset function.
Traditionally, the outset function is assumed to be homogeneous
\footnote{By definition, the homogeneous function meets the condition: $f(\lambda x)=\lambda^q f(x)$.}
function by default
in all preceding investigations described in the literature. As the applications show, this does not match the
situations in practice.

We first consider
the conditions which may lead to the nullification of $f_A$ or $f_S$ depending on the dimension of ${\bf x}$-space.
For the sake of simplicity, let us focus on the even dimension of space, $n=2$
(that is, $\vec{\bf x}\in \mathbb{R}^2$). In this case, (\ref{Inv-F-t-2-4-S}) and (\ref{Inv-F-t-2-4-A}) take the following forms
\begin{eqnarray}
\label{Inv-F-t-2-4-S2}
&&f_S(\vec{\bf x})\Big|_{\mathbb{R}^2} =
-\int_{-\infty}^{+\infty} (d\eta)\, \frac{\mathcal{P}}{\eta^2}\,
\int_{0}^{2\pi} d\varphi \, \mathcal{R}[f](\eta + \langle \vec{\bf n}_\varphi, \vec{\bf x}\rangle, \varphi)
\end{eqnarray}
and
\begin{eqnarray}
\label{Inv-F-t-2-4-A2}
&&f_A(\vec{\bf x})\Big|_{\mathbb{R}^2} =
-i\pi \int_{-\infty}^{+\infty} (d\eta) \,  \delta(\eta)\,  \frac{\partial}{\partial\eta} \,\int_{0}^{2\pi} d\varphi
\mathcal{R}[f](\eta + \langle \vec{\bf n}_\varphi, \vec{\bf x}\rangle, \varphi).
\end{eqnarray}

Before going further, it is necessary to notice that the direct Radon transforms can be treated as the curve-linear integrations
of first kind. Indeed, we can write the following \cite{Gelfand:1964}
\begin{eqnarray}
\label{DRT-1-1}
\mathcal{R}[f](\tau,\varphi)&=&
\int_{-\infty}^{+\infty} d^2 \vec{\bf x} \, f(\vec{\bf x})
\delta\left( \tau - \langle\vec{\bf n}_\varphi, \vec{\bf x}\rangle\right)
\nonumber\\
&=&
\int_{L(\tau, \varphi)}\, ds\, f(\tau \,\cos\varphi - s\,\sin\varphi, \tau\,\sin\varphi + s\, \cos\varphi),
\end{eqnarray}
where $L(\tau, \varphi)$ corresponds to the line integration and the corresponding rotation
of the coordinate system is given by
\begin{eqnarray}
\label{SYS-1}
&&\vec{\bf x}=(x_1, x_2) \quad \Rightarrow \quad \vec{\bf x}^\prime=(p, s),
\nonumber\\
&&
x_{1}(p, s; \varphi)=p\,\cos\varphi - s\,\sin\varphi, \quad
x_{2}(p, s; \varphi)=p\,\sin\varphi + s\,\cos\varphi, 
\nonumber\\
&&
dx_1\, dx_2 = |J|\,dp \,ds \quad \text{with}\quad |J|=1 .
\end{eqnarray}

We now concentrate on the most trivial and ideal case where the outset function, by definition, is homogeneous and {\it (a)} unbounded or
{\it (b)} the function support is limited and symmetric. In this case,
the Radon image possesses the symmetry property written as
\begin{eqnarray}
\label{RT-Sym-1}
\mathcal{R}[f](\tau,\varphi)= \mathcal{R}[f](-\tau,\varphi) \quad \text{but} \quad
\mathcal{R}[f](-\tau,\varphi)\equiv \mathcal{R}[f](\tau,\varphi+\pi),
\end{eqnarray}
or, equivalently,
\begin{eqnarray}
\label{RT-Sym-2}
\int_{L(\tau, \varphi)}\, ds\, f\big(x_1(\tau,s; \varphi),\, x_2(\tau,s; \varphi)\big) =
\int_{L(-\tau, \varphi)}\, ds\,  f\big(x_1(\tau,s; \varphi),\, x_2(\tau,s; \varphi)\big).
\end{eqnarray}
Hence, if the conditions of (\ref{RT-Sym-1}) and (\ref{RT-Sym-2}) take place and the angular integration covers the full interval,
we can readily see that
the contribution of $f_A\Big|_{\mathbb{R}^2}$ merely disappears from the consideration, see \cite{Anikin:2019oes}
for all deatals.

However, if we deal with the outset function which is inhomogeneous and/or the above-mentioned conditions
{\it (a)} and {\it (b)} have been broken, we have that (for definiteness,
$\vec{\bf x}\in L(\tau, \varphi)\in \{\Omega_I| x_{1,2}>0 \}$
while $\vec{\bf x}\in L(-\tau, \varphi)\in \{\Omega_{III}| x_{1,2}<0 \}$)
\begin{eqnarray}
\label{RT-Sym-3}
f_{I}(\vec{\bf x})\equiv f(\vec{\bf x})\Theta\big( \vec{\bf x}\in \Omega_I\big)=F_1(x_1, x_2), \quad
f_{III}(\vec{\bf x})\equiv f(\vec{\bf x})\Theta\big( \vec{\bf x}\in \Omega_{III}\big)=F_2(x_1, x_2)
\end{eqnarray}
and, as consequence,
\begin{eqnarray}
\label{RT-Sym-4}
\int_{L(\tau, \varphi)}\, ds\, f\big(x_1(\tau,s; \varphi),\, x_2(\tau,s; \varphi)\big) \not=
\int_{L(-\tau, \varphi)}\, ds\,  f\big(x_1(-\tau,s; \varphi),\, x_2(-\tau,s; \varphi)\big).
\end{eqnarray}
Hence, in the function $f_A$ of (\ref{Inv-F-t-2-4-A}), the $\eta$-integrand averaged over $\varphi$
can be presented as
\begin{eqnarray}
\label{Inv-F-t-2-4-A3}
&&f_A(\vec{\bf x})\Big|_{\mathbb{R}^2} \sim
\overline{\mathcal{R}}[f](\eta ; \vec{\bf x})
\stackrel{\text{def.}}{=}
\int_{0}^{2\pi} d\varphi\,
\mathcal{R}[f](\eta + \langle \vec{\bf n}_\varphi, \vec{\bf x}\rangle, \varphi)=
\nonumber\\
&&
\int_{0}^{\pi} d\varphi \, \Big\{  \mathcal{R}[f_{I}](\eta + \langle \vec{\bf n}_\varphi, \vec{\bf x}\rangle, \varphi) +
\mathcal{R}[f_{III}](-\eta + \langle \vec{\bf n}_\varphi, \vec{\bf x}\rangle, \varphi) \Big\} \not= 0
\end{eqnarray}
provided the condition of (\ref{RT-Sym-4}) and without the $\eta$-integration with the corresponding weight operation.
%
%Here, $\mathcal{R}[f]$ and $\widetilde{\mathcal{R}}[f]$ correspond to $F_1$- and $F_2$-functions, see (\ref{RT-Sym-3}).
%
Further, if we take into account the $\eta$-integration, we reach the same conclusion
regarding whether or not the contribution of $f_A$ nullifies. Indeed, we have
\begin{eqnarray}
\label{Inv-F-t-2-4-A4}
&&f_A(\vec{\bf x})\Big|_{\mathbb{R}^2} =
\int_{-\infty}^{+\infty} d\eta\, \delta(\eta)
\int_{0}^{2\pi} d\varphi\,
\mathcal{R}^\prime_\eta[f](\eta + \langle \vec{\bf n}_\varphi, \vec{\bf x}\rangle, \varphi)=
\nonumber\\
&&
\int_{-\infty}^{+\infty} d\eta\, \delta(\eta)
\int_{0}^{\pi} d\varphi \, \Big\{  \mathcal{R}^\prime_\eta[f_{I}](\eta + \langle \vec{\bf n}_\varphi, \vec{\bf x}\rangle, \varphi) +
\mathcal{R}^\prime_\eta[f_{III}](-\eta + \langle \vec{\bf n}_\varphi, \vec{\bf x}\rangle, \varphi) \Big\} =
\nonumber\\
&&
\int_{-\infty}^{+\infty} d\eta\, \delta(\eta)
\int_{0}^{\pi} d\varphi \, \Big\{  \mathcal{R}^\prime_\eta[f_{I}](\eta + \langle \vec{\bf n}_\varphi, \vec{\bf x}\rangle, \varphi) -
\mathcal{R}^\prime_\eta[f_{III}](\eta + \langle \vec{\bf n}_\varphi, \vec{\bf x}\rangle, \varphi) \Big\}
\not= 0
\end{eqnarray}
which is valid due to the inhomogeneity and
without any extra angular restrictions for the Radon transforms.
This is one of our key conclusions presented in the paper.

To conclude this section, it is worth to notice that in the odd dimension of space, for example in $\mathbb{R}^3$,
the contribution of $f_A$ begins to be as the main contribution while the term of $f_S$ disappears provided
the homogeneity and symmetry of outset function discussed above.

%%%%%%%%%%%%%%%%
\section{The angular restrictions and boundary}
\label{Ang-rest}
%%%%%%%%%%%%%%%%%%%%%%%%%%%%%%%%%%%%

In the section, we are focusing on the influence of the finite and restricted support, where the outset function has been determined,
on the boundary conditions and, in its turn, on the angular dependence of the Radon images.
In what follows the inhomogeneity of outset function is now not important.

We again adhere the two dimensional Euclidian $\vec{\bf x}$-space, $\mathbb{R}^2$, and we begin with the
Fourier transform written as
\begin{eqnarray}
\label{F-t-1-d-w1}
\mathcal{F}[f](\vec{\bf q})= \int_{-\infty}^{+\infty} d^2 \vec{\bf x} \,
e^{-i\langle\vec{\bf q},\vec{\bf x}\rangle} \,f(\vec{\bf x})
\Big\{
\int_{-\infty}^{+\infty} (dT)\, \delta\left( T - \langle\vec{\bf q},\vec{\bf x}\rangle  \right)
\Big\},
\end{eqnarray}
where the integral representation of unit has been inserted.

Let us rewrite the mentioned Fourier slice theorem, see (\ref{F-t-3-dir}), in the different form.
It reads
\begin{eqnarray}
\label{FST-1}
\mathcal{F}[f](p, q_2) =\int_{-\infty}^{+\infty}(dt) \,e^{-iq_2t}\,
\mathcal{R}[f](t,p),
\end{eqnarray}
where the Radon image is given by
\begin{eqnarray}
\label{FST-2}
\mathcal{R}[f](t, p)=
\int_{-\infty}^{+\infty} d^2 \vec{\bf x} \, f(\vec{\bf x})
\delta\left( t + px_1 - x_2 \right)
\end{eqnarray}
with
\begin{eqnarray}
\label{Par-FST-1}
t\stackrel{\text{def.}}{=}\frac{T}{q_2}, \quad p\stackrel{\text{def.}}{=}-\frac{q_1}{q_2}.
\end{eqnarray}
The slop parameter $p$ in the line paramtrization of (\ref{FST-2}) reflects, in the other words, the
angular dependence of both the Fourier and Radon transforms.
In the direct Radon transform (\ref{FST-2}), the delta-function argument gives the condition:
$x_2=px_1+t$ which is the standard parametrization of straightforward line.

We dwell on the case of restricted support of the outset function given by
\begin{eqnarray}
\label{OF-1}
f(\vec{\bf x}) \Longrightarrow f(\vec{\bf x}) \Theta\big( x_1\in \Omega_{\Box}\big),
\quad \Omega_\Box \stackrel{\text{def.}}{=} \big\{-1\leq x_1 \leq 1; -1\leq x_2 \leq 1 \big\},
\end{eqnarray}
where $\Theta$ stands for the corresponding theta-function and the components of $\vec{\bf x}$ are independent from each other.

\begin{figure}[t]
\centerline{\includegraphics[width=0.5\textwidth]{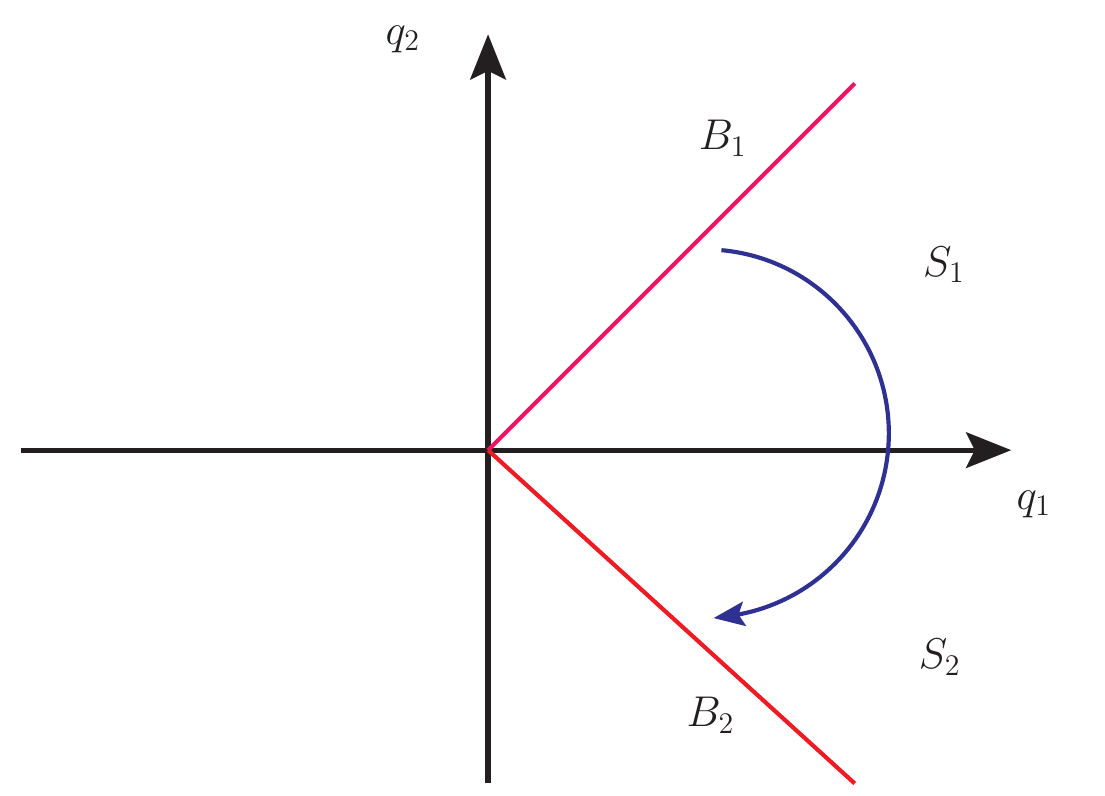}}
%\vspace{-0.5cm}
\caption{The restricted region in $\vec{\bf q}$-plane.
Notations: $B_1$ and $B_2$ correspond to the boundaries
given by $\{(1 + t)q_2=q_1 \,|\, q_{1,2}>0, t>0 \}$ and $\{ (1-t)q_2=q_1 \,|\, q_1>0, q_2<0, 1>t>0 \}$,
while the domains $S_1$ and $S_2$ correspond to $\{(1+t)q_2<q_1 \,|\, q_{1,2}>0, t>0 \}$
and $\{ (1-t)q_2<q_1 \,|\, q_1>0, q_2<0, 1>t>0 \}$.
}
\label{Fig-1}
\end{figure}

Due to the inequalities of (\ref{OF-1}) defining $\Omega_\Box$, we have
\begin{eqnarray}
\label{Inq-1}
-1 \leq px_1 + t\leq 1\quad \Rightarrow \quad
1 \stackrel{\text{\bf a}}{\geq} \frac{(t+1)q_2}{q_1} \geq x_1 \geq \frac{(t-1)q_2}{q_1}\stackrel{\text{\bf b}}{\geq} -1,
\end{eqnarray}
where the conditions $t>0$ and $\vec{\bf q}\in \{\Omega^I|\, q_{1,2}>0\}$ have been applied.
Hence, from (\ref{Inq-1}), the inequalities $\text{\bf a}$ and $\text{\bf b}$ readily give
the following conditions (see Fig.~\ref{Fig-1})
\begin{eqnarray}
\label{Inq-a-b-1}
  \begin{cases}
  (t+1)q_2 \leq q_1, &\text{for $\forall t$;}\\
  (1-t)q_2 \leq q_1, &\text{for $t<1$}
  \end{cases}
\end{eqnarray}
that restrict the variation domain for the variables $\vec{\bf q}$.

Focusing on the domain given by $\vec{\bf q}\in \{\Omega^{IV}|\, q_1>0, q_2<0 \}$,
we can similarly obtain that the inequalities $\text{\bf a}$ and $\text{\bf b}$ correspond to
(see Fig.~\ref{Fig-1})
\begin{eqnarray}
\label{Inq-a-b-2}
  \begin{cases}
  -(\tilde t-1) q_2 \leq q_1, &\text{for $\tilde t>1$, $\tilde t\equiv -t >0$}\\
  (1-t)q_2 \leq q_1, &\text{for $t<1$, $t>0$}.
  \end{cases}
\end{eqnarray}
The case of $\vec{\bf q}\in \{\Omega^{II}|\, q_1<0, q_2>0 \}\cup \{\Omega^{III}|\, q_1<0, q_2<0 \}$
corresponds to the restricted region that is mirrored left
to the restricted region presented in Fig.~\ref{Fig-1}. We omit this case in our discussion because
it does not give new information in the context of the (Radon) angular restrictions.

The special attention should be paid
for the boundary given by the equalities of (\ref{Inq-a-b-1}) and (\ref{Inq-a-b-2}),
{\it i.e.} described, for example, by $(t\pm1)q_2 = q_1$.
Indeed, these conditions reduce a number of independent Fourier (Radon) variables to one.
At the same time, the boundary of outset function has been still formed by two independent variables.

In the mathematical literature, the following theorem on the correspondence of numbers of independent integration variables
is well-known:
{\it if the outset function $f$ is a function of $N$ independent variables then the Radon
image $\mathcal{R}[f]$ depends on $N$ independent variables too.
The Radon transforms are the bijections and exist on
the space of $\mathbb{R}^1\times \mathrm{S}^{N-1}$} \cite{Deans}.
It means that to avoid the discrepancy between the boundary transformations of outset function and the
Fourier (Radon) image, we have to exclude the boundary corresponding to the Fourier (Radon) image
from the consideration. That is, a crossing via the boundary $B_1$ and $B_2$ from the domains $S_1$
and $S_2$ is forbidden, see Fig.~\ref{Fig-1}.
As a result, the angular dependence of Fourier (Radon) transform
receives the definite limits given by the maximal interval $(-\pi/2; \, \pi/2)$.

Thus, we have shown that the restricted support of outset function produces the angular restriction
for the Radon images. In its turn, it leads to the essential contributions of both terms $f_S$ and $f_A$
to the inversion of Radon transform, see (\ref{Inv-F-t-2-4-S}) and (\ref{Inv-F-t-2-4-A}).

%%%%%%%%%%%%%%%%%%%%%%%%%%%%% FIGURE %%%%%%%%%%%%%%%%%%%%%%%%%%%%%%%%
%%%%%%%%%%%%%%%%%%%%%%%%%%%%%%%%%%%%%%%%%%%%%%%%%%%%%%%%%%%%%%%%%%%%%%%

%%%%%%%%%%%%%%%%
\section{The surface term and boundary}
\label{Serf-BC}
%%%%%%%%%%%%%%%%%%%%%%%%%%%%%%%%%%%

We are now in a position to show that the exclusion of boundary from the support domain, where the outset function has been defined,
is dictated by the requirement of the surface term absence, or vise versa.
For the simplicity, but without loosing the generality, we again work in
$\mathbb{R}^2$-space.

In the generalized function theory, the surface terms can be appeared, in particular, in integration by part if the integrand involves the
derivative of generalized delta-function. Indeed, let us consider the simplest example given by
the following integration:
\begin{eqnarray}
\label{Delta-1}
\int_{\mathcal{D}} (dx) \, \varphi(x)\, \partial_x \delta(x) =
\varphi(x)\, \delta(x) \Big|_{\text{B.}} -
\int_{\mathcal{D}} (dx) \, \delta(x)\, \partial_x \varphi(x),
\end{eqnarray}
where $\varphi\in \mathcal{D}$ is a finite function by definition; $``\text{B.}"$ denotes the integration limits (or the boundary).
Concentrating on the surface (first) term of {\it r.h.s.} of (\ref{Delta-1}),
one can see that the point $x=0$ corresponding to the delta-function argument is, as usual, out
of the integration limits (boundary). Therefore, the surface term in (\ref{Delta-1}) does not contribute at all
provided the function $\varphi$ belongs to the set of finite functions, $\mathcal{D}$.
In this context, we want to study the Radon transformation of outset function which is also
involving the delta-function as a part of integrand.

In both $f_S$ and $f_A$ of (\ref{Inv-F-t-2-4}), we deal explicitly or implicitly with the derivatives over the radial Radon
component, $\partial_\tau$, that act on the Radon image.
Indeed, we have the following typical combination:
\begin{eqnarray}
\label{ST-1}
&&\int d\mu(\tau, \varphi)\frac{\partial}{\partial\tau}
\mathcal{R}[f]\left( \tau + \langle\vec{\bf n}_\varphi, \vec{\bf x} \rangle, \varphi  \right)=
\nonumber\\
&&
-\int d\mu(\tau, \varphi) \int_{\Omega} d^2 \vec{\bf y}\, f(\vec{\bf y}) \,
\langle\vec{\bf n}_\varphi, \vec{\nabla}_y \rangle
\delta\left( \tau + \langle\vec{\bf n}_\varphi, \vec{\bf x}- \vec{\bf y} \rangle \right),
\end{eqnarray}
where the integration measure is symbolically presented through $d\mu(\tau, \varphi)$ and
$\partial_\tau \Rightarrow - \langle\vec{\bf n}_\varphi, \vec{\nabla}_y \rangle$.
Next, we make a replacement: $\vec{\bf x}-\vec{\bf y}=\vec{\bf z}$ and, then, we integrate by part leading to
\begin{eqnarray}
\label{ST-2}
&&- \int d\mu(\varphi) \int_{\Omega} d^2 \vec{\bf z}\, \big[ \langle\vec{\bf n}_\varphi, \vec{\nabla}_z \rangle
f(\vec{\bf x}-\vec{\bf z}) \big] \,
\delta\left( \tau + \langle\vec{\bf n}_\varphi, \vec{\bf z} \rangle \right)+
\nonumber\\
&&\int d\mu(\varphi; \cos\varphi) \int_{\Omega_2} dz_2
f(\vec{\bf x}-\vec{\bf z})\,
\delta\left( \tau + \langle\vec{\bf n}_\varphi, \vec{\bf z} \rangle \right) \Big|_{\Omega_1} +
\nonumber\\
&&\int d\mu(\varphi; \sin\varphi) \int_{\Omega_1} dz_1
f(\vec{\bf x}-\vec{\bf z}) \,
\delta\left( \tau + \langle\vec{\bf n}_\varphi, \vec{\bf z} \rangle \right) \Big|_{\Omega_2},
\end{eqnarray}
where $\Omega_{1,2}$ denote the boundary with respect to $z_1$ and $z_2$.
For brevity, the surface terms of (\ref{ST-2}) can be presented in the simplified forms as
\begin{eqnarray}
\label{ST-3}
\int d\mu_F(\varphi)
F_{\,\vec{\bf x}}(\vec{\bf z})\,
\delta\left( \tau + \langle\vec{\bf n}_\varphi, \vec{\bf z} \rangle \right) \Big|_{\Omega_F},
\end{eqnarray}
where $F_{\,\vec{\bf x}}(\vec{\bf z})\equiv f(\vec{\bf x}-\vec{\bf z})$ and, for our aims,
the integrations over the components of
$\vec{\bf z}$ have been omitted as they become irrelevant.

The radial $\tau$ dependence of Radon image is necessarily restricted if the outset function has been
well-localized in ${\vec{\bf x}}$-space \cite{Deans}. This is a reason of the surface $\tau$-term absence
in the {\it l.h.s.} of (\ref{ST-1}). The implemented replacement: $\partial_\tau \Rightarrow \vec{\nabla}_y$
has to keep this mentioned property to be valid too.
Therefore, it leads to the inference that the surface terms of ({\ref{ST-3}}) have to be disappeared as well.

\begin{figure}[t]
\centerline{\includegraphics[width=0.5\textwidth]{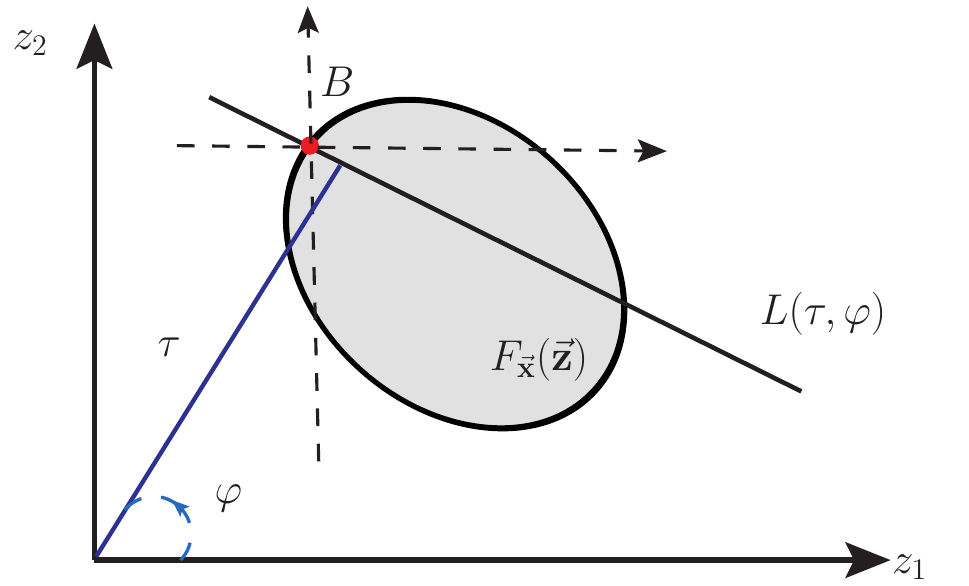}}
%\vspace{-0.5cm}
\caption{The exclusion of boundary that leads to the surface term absence.
The red point $B$ belongs to both the domain boundary
and the line defined the Radon transform.
}
\label{Fig-2}
\end{figure}

Next, in (\ref{ST-3}), the argument of delta-function defines the parametrization of straightforward line.
From the practical viewpoint,
we suppose that the outset-like function $F_{\,\vec{\bf x}}(\vec{\bf z})$ has been restricted by some
domain with the non-zero support, see Fig.~\ref{Fig-2}.
The straightforward line defined by $\delta$-argument is always crossing the boundary of the domain where the outset-like function
has been determined. To avoid the existence of surface terms, which is artificial for our case,
we have to exclude the boundary of domain.
In other words, the non-zero support of outset function has been formed by the open domain.
As mentioned, the exclusion of boundary agrees with the theorem on the correspondence of
numbers of independent variables for the outset function and its Radon image.

%%%%%%%%%%%%%%%%%%%%%%%%%%%%% FIGURE %%%%%%%%%%%%%%%%%%%%%%%%%%%%%%%%
%%%%%%%%%%%%%%%%%%%%%%%%%%%%%%%%%%%%%%%%%%%%%%%%%%%%%%%%%%%%%%%%%%%%%%%

%%%%%%%%%%%%%%%%
\section{Conclusions}
\label{Cons}
%%%%%%%%%%%%%%%%%%%%%%%%%%%%%%%%%%%%

To conclude, we have obtained the universal expression of Radon inversion
which can be used for both even and odd dimensions simultaneously.
This has been achieved with the help of the suitable
regularization in the frame of the generalized function and
without the use of Courant-Hilbert's identities.

We have demonstrated that, in the universal representation of inverse Radon transform,
there are two
essential terms, $f_S$ and $f_A$, which are contributing
even in the full angular region of variation, {\it i.e.} in the interval $(0, 2\pi)$,
due to the inhomogeneity of outset function.

In the paper, the restrictions on the Radon angular dependence have been derived.
These restrictions ensure the contributions of $f_S$ and $f_A$
independently from the homogeneous properties of outset function.
It has been shown that the restricted angular dependence of Radon transforms
is a consequence of the finite support of outset functions localized in the space domain.
Also, we have proved that the angular restrictions
are closely related to the
the boundary exclusion of outset function and, therefore, of the Radon images.

%%%%%%%%%%%%%%
\section*{Acknowledgements}
%%%%%%%%%%%%%%
We thank M.V.~Kompaniets, V.A.~Osipov, L.~Szymanowski and N.A.~Tiurin for useful and illuminating discussions.
The special thanks go to the colleagues from the South China Normal University (Guangzhou) and
the Chinese University of Hong Kong (Shenzhen) for the useful discussions and a very warm hospitality.

\section*{Conflict of interest}

The author declares no conflict of interest.

\section*{Data Availability Statement}
This manuscript has no associated data or the data will not be deposited.

%%%%%%%%%%%%%%%%%%%%%%%%%%%

\end{document}